\documentclass[10pt,a4paper]{article}
\usepackage{cite}
\usepackage{amscd}
\usepackage[latin1]{inputenc}
\usepackage{amsmath}
\usepackage{amsfonts}
\usepackage{amssymb}
\usepackage{color}
\usepackage{float}
\usepackage{amsthm}
\usepackage{graphicx}
\usepackage{listings}
\usepackage{hyperref}
\usepackage{mathtools}

\newtheorem{theorem}{Theorem}

\newtheorem{proposition}[theorem]{Proposition}
\newtheorem{remark}[theorem]{Remark}

\def\m{\mathbb}
\def\mc{\mathcal}

\begin{document}
\date{}
\title{A surface with canonical map of degree $24$}
\author{Carlos Rito}
\maketitle

\begin{abstract}

We construct a complex algebraic surface with geometric genus $p_g=3$, irregularity $q=0$, self-intersection of the canonical divisor $K^2=24$ and canonical map of degree $24$ onto $\m P^2$.


\end{abstract}

\section{Introduction}

Let $S$ be a smooth minimal surface of general type with geometric genus $p_g\geq 3$.
Denote by ${\phi:S\dashrightarrow\m P^{p_g-1}}$ the canonical map and let $d:=\deg(\phi).$
The following Beauville's result is well-known.
\begin{theorem}[\cite{Be}]
If the canonical image $\Sigma:=\phi(S)$ is a surface, then either:
\begin{description}
\item{{\rm (i)}} $p_g(\Sigma)=0,$ or
\item{{\rm (ii)}} $\Sigma$ is a canonical surface $($in particular $p_g(\Sigma)=p_g(S))$.
\end{description}
Moreover, in case {\rm (i)} $d\leq 36$ and in case {\rm (ii)} $d\leq 9.$
\end{theorem}

Beauville has also constructed families of examples with $\chi(\mathcal O_S)$ arbitrarily
large for $d=2, 4, 6, 8$ and $p_g(\Sigma)=0.$ Despite being a classical problem, for $d>8$ the
number of known examples drops drastically. Tan's example \cite[\S 5]{Ta} with
$d=9$ and Persson's example \cite{Pe} with $d=16,$ $q=0$ are well known.
Du and Gao \cite{DG} show that if the canonical map is an abelian cover
of $\m P^2,$ then these are the only possibilities for $d>8.$
More recently the author has given examples with $d=12$
\cite{Ri2} and $d=16,$ $q=2$ \cite{Ri3}.

In this paper we construct a surface $S$ with $p_g=3,$ $q=0$ and $d=24,$ obtained as a
$\m Z_2^4$-covering of $\m P^2.$
The canonical map of $S$ factors through a $\m Z_2^2$-covering of a surface with
$p_g=3,$ $q=0$ and $K^2=6$ having $24$ nodes, which in turn is a double covering of a Kummer surface.


\bigskip
\noindent{\bf Notation}

We work over the complex numbers. All varieties are assumed to be projective algebraic.
A $(-n)$-curve on a surface is a curve isomorphic to $\m P^1$ with self-intersection $-n.$
Linear equivalence of divisors is denoted by $\equiv.$
The rest of the notation is standard in Algebraic Geometry.\\

\bigskip
\noindent{\bf Acknowledgements}

The author would like to thank an anonymous referee for suggestions to improve the exposition of the paper.

This research was partially supported by FCT (Portugal) under the project PTDC/MAT-GEO/0675/2012, the fellowship SFRH/BPD/111131/2015
and by CMUP (UID/MAT/00144/2013), which is funded by FCT with national (MEC) and European structural funds through the programs FEDER, under the partnership agreement PT2020.

\section{$\m Z_2^n$-coverings}\label{coverings}

The following is taken from \cite{Ca}, the standard reference is \cite{Pa}.

\begin{proposition}
A normal finite $G\cong\m Z_2^r$-covering $Y\rightarrow X$ of a smooth variety $X$ is completely
determined by the datum of
\begin{enumerate}
\item reduced effective divisors $D_{\sigma},$ $\forall\sigma\in G,$ with no common components;
\item divisor classes $L_1,\ldots,L_r,$ for $\chi_1,\ldots,\chi_r$ a basis of the dual group of characters
$G^{\vee},$ such that $$2L_i\equiv\sum_{\chi_i(\sigma)=-1}D_{\sigma}.$$ 
\end{enumerate}
Conversely, given 1. and 2., one obtains a normal scheme $Y$ with a finite
$G\cong\m Z_2^r$-covering $Y\rightarrow X,$ with branch curves the divisors $D_{\sigma}.$
\end{proposition}
The covering $Y\rightarrow X$ is embedded in the total space of the direct sum of the line bundles whose
sheaves of sections are the $\mc O_X(L_i),$ and is there defined by equations
$$u_{\chi_i}u_{\chi_j}=u_{\chi_i\chi_j}\prod_{\chi_i(\sigma)=\chi_j(\sigma)=-1}x_{\sigma},$$
where $x_{\sigma}$ is a section such that ${\rm div}(x_{\sigma})=D_{\sigma}.$

The scheme $Y$ can be seen as the normalization of the Galois covering given by the equations
$$u_{\chi_i}^2=\prod_{\chi_i(\sigma)=-1}x_{\sigma},$$
and $Y$ is irreducible if $\{\sigma|D_{\sigma}>0\}$ generates $G.$\\

For a covering $\pi:Y\rightarrow X$ with ramification divisor $R,$ the Hurwitz formula $K_Y=\pi^*(K_X)+R$ holds.
Let us describe the canonical system for the case where $\pi$ is a $\mathbb Z_2^2$-covering with smooth branch divisor.
We have branch curves $D_1,D_2,D_3$ and relations $2L_i\equiv D_j+D_k,$ for all permutations $(i,j,k)$ of $\{1,2,3\}$.
The covering $\pi$ factors as $$\phi:Y\rightarrow W_i,\ \ \varphi:W_i\rightarrow X,$$
where $\varphi$ is the double covering corresponding to $L_i.$
Let $R_i$ be the ramification divisor of $\phi.$ One has
$$K_Y\equiv\phi^*(K_{W_i})+R_i\ \ \ {\rm and}\ \ \ K_{W_i}\equiv\varphi^*(K_X+L_i),$$
which gives
$$K_Y\equiv\pi^*(K_X+L_i)+\frac{1}{2}\pi^*(D_i),\ \ i=1,2,3.$$

Finally we notice that taking the quotient by a subgroup $H$ of the
Galois group of the covering corresponds to considering the subalgebra generated by the line bundles
$L_{\chi}^{-1},$ where $\chi$ ranges over the characters orthogonal to $H.$

\section{The construction}\label{section3}

We show in the Appendix the existence of reduced plane curves $C_6$ of degree $6$ and $C_7$ of degree $7$
through points $p_0,\ldots,p_5$ such that:
\begin{description}
\item{$\cdot$} $C_7$ has a triple point at $p_0$ and tacnodes at $p_1,\ldots,p_5;$
\item{$\cdot$} $C_6$ is smooth at $p_5,$ has a node at $p_0$ and tacnodes at $p_1,\ldots,p_4;$
\item{$\cdot$} the branches of the tacnode of $C_j$ at $p_i$ are tangent to the line $T_i$ through $p_0,p_i,$ $j=1,2,$ $i=1,\ldots,4;$
\item{$\cdot$} the branches of the tacnode of $C_7$ at $p_5$ are tangent to $C_6;$
\item{$\cdot$} the singularities of $C_6+C_7$ are resolved via one blow-up at $p_0$ and two blow-ups at each of $p_1,\ldots,p_5.$
\end{description}

\noindent{\bf Step 1} (Construction)\\
Consider the map $$\mu:X\longrightarrow\m P^2$$ which resolves the singularities of the curve $C_7.$
Then $\mu$ is given by blow-ups at $$p_0,p_1,p_1',\ldots,p_5,p_5',$$ where $p_i'$ is infinitely near to $p_i.$
Let $E_0,E_1,E_1',\ldots,E_5,E_5'$ be the corresponding exceptional divisors (with self-intersection $-1$).

Let $x, y, z, w$ be generators of the group $\m Z_2^4$ and $$\psi:Y\longrightarrow X$$ be the $\m Z_2^4$-covering defined by
$$D_x:=\widetilde T_1-E_0-2E_1',$$
$$D_y:=\widetilde T_2-E_0-2E_2',$$
$$D_z:=\widetilde C_6-2E_0-\sum_1^4(2E_i+2E_i')-2E_5',$$
$$D_w:=\widetilde C_7+\widetilde T_4-4E_0-\sum_1^3(2E_i+2E_i')-(2E_4+4E_4')-(2E_5+2E_5'),$$
$$D_{xy}:=\widetilde T_3-E_0-2E_3',$$
$$D_{xz}:=\cdots:=D_{zw}:=0,$$
where the notation $\widetilde{\cdot}$ stands for the total transform $\mu^*(\cdot).$

We note that each of the divisors $D_x,$ $D_y,$ $D_{xy}$ and $\widetilde T_4-E_0-2E_4'$ (contained in $D_w$)
is a disjoint union of two $(-2)$-curves.

For $i,j,k,l\in\{-1,1\},$ let $\chi_{ijkl}$ denote the character which takes the value $i,j,k,l$ on $x,y,z,w,$ respectively.
There exist divisors $L_{ijkl}$ such that
\begin{equation}\label{eq1}
2L_{ijkl}\equiv\sum_{\chi_{ijkl}(\sigma)=-1}D_{\sigma},
\end{equation}
thus the covering $\psi$ is well defined.
Since there is no 2-torsion in the Picard group of $X,$ then $\psi$ is uniquely determined.
The surface $Y$ is smooth because the curves $D_x,\ldots,D_{xy}$ are smooth and disjoint.
Division of the equations (\ref{eq1}) by 2 gives that the $L_{ijkl}$ are according to the following table.
For instance $L_{-1111}\equiv \widetilde T-E_0-E_1'-E_3'.$
{\small
$$
\bordermatrix{~ & \widetilde T & E_0 & E_1 & E_1' & E_2 & E_2' & E_3 & E_3' & E_4 & E_4' & E_5 & E_5' \cr
L_{-1111} & 1 & -1 & 0 & -1 & 0 & 0 & 0 & -1 & 0 & 0 & 0 & 0 \cr
L_{1-111} & 1 & -1 & 0 & 0 & 0 & -1 & 0 & -1 & 0 & 0 & 0 & 0 \cr
L_{-1-111} & 1 & -1 & 0 & -1 & 0 & -1 & 0 & 0 & 0 & 0 & 0 & 0 \cr
L_{11-11} & 3 & -1 & -1 & -1 & -1 & -1 & -1 & -1 & -1 & -1 & 0 & -1 \cr
L_{-11-11} & 4 & -2 & -1 & -2 & -1 & -1 & -1 & -2 & -1 & -1 & 0 & -1 \cr
L_{1-1-11} & 4 & -2 & -1 & -1 & -1 & -2 & -1 & -2 & -1 & -1 & 0 & -1 \cr
L_{-1-1-11} & 4 & -2 & -1 & -2 & -1 & -2 & -1 & -1 & -1 & -1 & 0 & -1 \cr
L_{111-1} & 4 & -2 & -1 & -1 & -1 & -1 & -1 & -1 & -1 & -2 & -1 & -1 \cr
L_{-111-1} & 5 & -3 & -1 & -2 & -1 & -1 & -1 & -2 & -1 & -2 & -1 & -1 \cr
L_{1-11-1} & 5 & -3 & -1 & -1 & -1 & -2 & -1 & -2 & -1 & -2 & -1 & -1 \cr
L_{-1-11-1} & 5 & -3 & -1 & -2 & -1 & -2 & -1 & -1 & -1 & -2 & -1 & -1 \cr
L_{11-1-1} & 7 & -3 & -2 & -2 & -2 & -2 & -2 & -2 & -2 & -3 & -1 & -2 \cr
L_{-11-1-1} & 8 & -4 & -2 & -3 & -2 & -2 & -2 & -3 & -2 & -3 & -1 & -2 \cr
L_{1-1-1-1} & 8 & -4 & -2 & -2 & -2 & -3 & -2 & -3 & -2 & -3 & -1 & -2 \cr
L_{-1-1-1-1} & 8 & -4 & -2 & -3 & -2 & -3 & -2 & -2 & -2 & -3 & -1 & -2 \cr}
$$
}

\noindent{\bf Step 2} (Invariants)\\
Since $$K_X\equiv -3\widetilde T+E_0+\sum_1^5(E_i+E_i'),$$ then
$$\chi(\mathcal O_Y)=16\chi(\mathcal O_X)+\frac{1}{2}\sum \left(L_{ijkl}^2+K_XL_{ijkl}\right)=$$ $$=16-1-1-1+0-1-1-1+0-1-1-1+0-1-1-1=4.$$

For the computation of $$p_g(Y)=p_g(X)+\sum h^0(X,\mathcal O_X(K_X+L_{ijkl})),$$
let
$$\mathcal T_1:=\left(\widetilde T_4-E_0-2E_4'+E_5-E_5'\right),$$
$$\mathcal T_2:=\left(\widetilde T_2+\widetilde T_3+\widetilde T_4-3E_0-\sum_2^4 2E_i'+E_5-E_5'\right),$$
$$\mathcal L_1:=\left|3\widetilde T-E_0-\sum_1^3(E_i+E_i')-E_4-E_5\right|$$
and
$$\mathcal L_2:=\left|2\widetilde T-(E_1+E_1')-E_2-E_3-E_4-E_5\right|.$$
Each of $\mathcal T_1,$ $\mathcal T_2$ is a disjoint union of $(-2)$-curves intersecting negatively $K_X+L_{11-1-1},$
$K_X+L_{1-1-1-1},$ respectively, thus we have
$$|K_X+L_{11-1-1}|=\mathcal T_1+\mathcal L_1$$
and
$$|K_X+L_{1-1-1-1}|=\mathcal T_2+\mathcal L_2.$$
We show in the Appendix that $\mathcal L_1$ has only one element and $\mathcal L_2$ is empty.
Hence
$$h^0(X,\mathcal O_X(K_X+L_{11-1-1}))=1$$ and
$$h^0(X,\mathcal O_X(K_X+L_{1-1-1-1}))=0.$$
Analogously 
$$h^0(X,\mathcal O_X(K_X+L_{-11-1-1}))=h^0(X,\mathcal O_X(K_X+L_{-1-1-1-1}))=0.$$

It is easy to see that
$$h^0(X,\mathcal O_X(K_X+L_{11-11}))=h^0(X,\mathcal O_X(K_X+L_{111-1}))=1$$
and $$h^0(X,\mathcal O_X(K_X+L_{ijkl}))=0$$ for the remaining cases.
We conclude that $$p_g(Y)=0+1+1+1=3.$$

Now we compute the self-intersection of the canonical divisor for the minimal model $S$ of $Y.$
The divisor $$\xi_1:=\frac{1}{2}\psi^*\left(\sum_1^3\left(\widetilde T_i-E_0-2E_i'\right)\right)$$
is a disjoint union of $8\times 6=48$ $(-1)$-curves
and the divisor $$\xi_2:=\frac{1}{2}\psi^*\left(\widetilde T_4-E_0-2E_4'+E_5-E_5'\right)$$
is a disjoint union of $8\times 3=24$ $(-1)$-curves.

The covering $\psi$ factors through the double covering $\varphi:W\rightarrow X$ with branch locus $D_z+D_w.$
We have $K_W\equiv\varphi^*(K_X+L_{11-1-1}),$ hence the Hurwitz formula gives
$$K_Y\equiv\xi_1+\psi^*(K_X+L_{11-1-1}).$$
Thus one of the canonical curves of $Y$ is $$\xi_1+2\xi_2+\psi^*(\mathcal C),$$
where $\mathcal C$ is the unique element in the linear system $\mathcal L_1$ defined above.
From $\xi_1\xi_2=\xi_1\psi^*(\mathcal C)=\psi^*(\mathcal C)^2=0$ and $\xi_2\psi^*(\mathcal C)=24,$ we get $K_Y^2=-48.$
We show in the Appendix that the curve $\mathcal C$ is irreducible, therefore $\psi^*(\mathcal C)$ is nef and then $K_S^2=24.$\\

\noindent{\bf Step 3} (The canonical map)\\
The divisors $$D_z,\ D_w,\ D_{zw}$$ define a $\mathbb Z_2^2$-covering
$$\rho:U\rightarrow X.$$
We have
$$\chi(\mathcal O_U)=4\chi(\mathcal O_X)+\frac{1}{2}\sum \left(L_{11kl}^2+K_XL_{11kl}\right)=4+0+0+0=4$$
and
$$p_g(U)=p_g(X)+\sum h^0(X,\mathcal O_X(K_X+L_{11kl}))=0+1+1+1=3.$$
The surface $U$ is the quotient of $Y$ by the subgroup $H$ generated by $x, y.$
The group $H$ acts on the minimal model $S$ of $Y$ with only isolated fixed points,
so $S/H$ is the canonical model $\bar U$ of $U$ and then $$K_{\bar U}^2=6.$$

Finally we want to show that the canonical map of $U$ is of degree $6$ onto $\m P^2.$
It suffices to verify that the canonical system has no base component nor base points.
The canonical system of $U$ is generated by the divisors
\begin{align*} 
K_1 &:=\frac{1}{2}\rho^*(D_z)+\rho^*(K_X+L_{111-1}),\\
K_2 &:=\frac{1}{2}\rho^*(D_w)+\rho^*(K_X+L_{11-11}),\\
K_3 &:=\frac{1}{2}\rho^*(D_{zw})+\rho^*(K_X+L_{11-1-1}).
\end{align*}
Denote by $\vartheta_1,\ldots,\vartheta_4$ the four $(-1)$-curves
$$\frac{1}{2}\rho^*(\widetilde T_4-E_0-2E_4')$$
and by $\vartheta_5,\vartheta_6$ the two $(-1)$-curves
$$\frac{1}{2}\rho^*(E_5-E_5').$$
Let $$\pi:U\rightarrow U'$$ be the contraction to the minimal model and
$q_1,\ldots,q_6\in U'$ be the points obtained by contraction of $\vartheta_1,\ldots,\vartheta_6.$
If $\kappa$ is an effective canonical divisor of $U',$ then $$H:=\pi^*(\kappa)+\vartheta_1+\cdots+\vartheta_6$$ is a canonical curve of $U.$
So, the multiplicity of a curve $\vartheta_i$ in $H$ is 1 if and only if the curve $\kappa$ does not contain the point $q_i.$

Since the multiplicity of $\vartheta_5+\vartheta_6$ in $K_1$ is $1$, the points $q_5,q_6$ are not
base points of the canonical system of $U'$.
The multiplicity of $\vartheta_1+\cdots+\vartheta_4$ in $K_2$ is $1$, so also the points $q_1,\ldots,q_4$
are not base points of the canonical system of $U'$.
Now to conclude the non-existence of other base points, it suffices to show that the plane curves
$$\mu\circ\rho(K_i), \ i=1,2,3,$$ have common intersection $\{p_0,p_1,\ldots,p_5\}$ and their
singularities are no worse than stated. This is done in the Appendix.
Here we just note that these curves are
$$T_4+C_6,\ \ C_7,\ \ T_4+C_3,$$
where $C_3$ is the plane cubic corresponding to the unique element in the linear system $\mathcal L_1,$ defined in Step 2 above.\\

\noindent{\bf Step 4} (Conclusion)\\
The $\m Z_2^4$-covering $\psi:Y\rightarrow X$ factors as
$$Y\xrightarrow{\ 4:1\ } U\xrightarrow{\ 4:1\ } X.$$ 
Since $p_g(Y)=p_g(U)=3$ and the canonical map of $U$ is of degree $6,$ then the canonical map of
$Y$ is of degree $24.$

\begin{remark}
Consider the intermediate double covering $\epsilon:Q\rightarrow X$ of $\rho$ with branch locus $D_z.$
Then $Q$ is a Kummer surface: each divisor
\mbox{$\epsilon^*\left(\widetilde T-E_0-2E_i'\right)$}
is a disjoint union of four $(-2)$-curves.
The surface $U$ contains $24$ disjoint $(-2)$-curves $A_1,\ldots,A_{24},$
the pullback of $\sum_1^3\epsilon^*\left(\widetilde T_i-E_0-2E_i'\right),$
such that the covering $Y\rightarrow U$ is a $\m Z_2^2$-Galois covering ramified over the divisors
$$A_1+\cdots+A_8, \ A_9+\cdots+A_{16}, \ A_{17}+\cdots+A_{24}.$$
\end{remark}

\appendix
\section*{Appendix}\label{appendix}
The following code is implemented with the Computational Algebra System Magma \cite{BCP},
version V2.21-8.\\

First we compute the curves $C_6$ and $C_7$ referred in Section \ref{section3}. We choose the points $p_0,\ldots,p_5$
with a symmetry axis and compute the curves using the Magma function $LinSys$ given in \cite{Ri1}.
\begin{verbatim}
A<x,y>:=AffineSpace(Rationals(),2);
P:=[A![0,0],A![2,2],A![-2,2],A![3,1],A![-3,1],A![0,5]];
M1:=[[2],[2,2],[2,2],[2,2],[2,2],[1,1]];
M2:=[[3],[2,2],[2,2],[2,2],[2,2],[2,2]];
T:=[[],[[1,1]],[[-1,1]],[[3,1]],[[-3,1]],[[1,0]]];
J6:=LinSys(LinearSystem(A,6),P,M1,T);
J7:=LinSys(LinearSystem(A,7),P,M2,T);
C6:=Curve(A,Sections(J6)[1]);
C7:=Curve(A,Sections(J7)[1]);
\end{verbatim}
\noindent We consider the projective closure of the curves and verify that they are irreducible and the singularities are exactly as stated.
\begin{verbatim}
P2<x,y,z>:=ProjectiveClosure(A);
C6:=ProjectiveClosure(C6);
C7:=ProjectiveClosure(C7);
IsAbsolutelyIrreducible(C6);
IsAbsolutelyIrreducible(C7);
SingularPoints(C6 join C7);
HasSingularPointsOverExtension(C6 join C7);
[ResolutionGraph(C6,P[i]):i in [1..#P-1]];
[ResolutionGraph(C7,P[i]):i in [1..#P]];
[ResolutionGraph(C6 join C7,P[i]):i in [1..#P]];
\end{verbatim}
\noindent
To clarify the situation at the origin, we use:
\begin{verbatim}
d:=DefiningEquation(TangentCone(C7,A![0,0]));
d eq y*(x^2 + 40585383/1587545*y^2);
\end{verbatim}
thus the singularity is ordinary.\\

\noindent The defining polynomials of $C_6$ and $C_7$ are
{\small
\begin{verbatim}
289*x^6+754326*x^4*y^2+2610657*x^2*y^4+1906344*y^6-2013848*x^4*y*z
-17946576*x^2*y^3*z-22212504*y^5*z+1336400*x^4*z^2
+35856160*x^2*y^2*z^2+89326224*y^4*z^2-22270208*x^2*y*z^3
-146421504*y^3*z^3+295936*x^2*z^4+84049920*y^2*z^4
\end{verbatim}
}
\noindent and
{\small
\begin{verbatim}
8683464*x^6*y-494984955*x^4*y^3-1064093674*x^2*y^5-558251235*y^7
-11358312*x^6*z+1253331746*x^4*y^2*z+8340957732*x^2*y^4*z
+7286240034*y^6*z-920312219*x^4*y*z^2-17394911410*x^2*y^3*z^2
-32292289971*y^5*z^2+179839940*x^4*z^3+11716330200*x^2*y^2*z^3
+55580514660*y^4*z^3-1270036000*x^2*y*z^4-32468306400*y^3*z^4
\end{verbatim}
}

Now we show that the linear system $\mathcal L_1,$ defined in Step 2 above, has exactly one element.
Let $L_1$ be the corresponding linear system of plane cubics.
By parameter counting, $\dim(L_1)\geq 0.$ If $\dim(L_1)\geq 1,$
then one of its curves contains the line $T_3,$ because
$$\left(\widetilde T_3-E_0-E_3-E_3'\right)\left(3\widetilde T-E_0-\sum_1^3(E_i+E_i')-E_4-E_5\right)=0.$$
The other component of this curve is a conic, but one can verify that the conic through $p_4$ tangent
to the lines $T_1,T_2$ at $p_1,p_2,$ which is given by the equation
\begin{displaymath}
x^2 - 9y^2 + 32y - 32=0,
\end{displaymath}
does not contain the point $p_5.$
We compute the unique plane cubic $C_3$ in $L_1$ and show that it is irreducible:
\begin{verbatim}
M:=[[1],[1,1],[1,1],[1,1],[1,0],[1,0]];
J3:=LinSys(LinearSystem(A,3),P,M,T);
#Sections(J3) eq 1;
C3:=ProjectiveClosure(Curve(A,Sections(J3)[1]));
IsAbsolutelyIrreducible(C3);
\end{verbatim}
\noindent The defining polynomial of $C_3$ is
{\small
\begin{verbatim}
17*x^3-924*x^2*y-153*x*y^2-996*y^3+1164*x^2*z
+544*x*y*z+6516*y^2*z-544*x*z^2-7680*y*z^2
\end{verbatim}
}

To conclude that the linear system $\mathcal L_2,$ defined in Step 2, is empty,
it suffices to note that the conic $C$ through $p_1,\ldots,p_5$ 
is not tangent to the line $T_1$ at the point $p_1.$
An equation for $C$ is
\begin{displaymath}
-12x^2 + 11y^2 - 93y + 190 = 0.
\end{displaymath}

Finally we verify that the curves
$$T_4+C_6,\ \ C_7,\ \ T_4+C_3,$$ referred in the end of Section \ref{section3},
have intersection $\{p_0,p_1,\ldots,p_5\}:$
\begin{verbatim}
T4:=Curve(P2,x+3*y);
PointsOverSplittingField((T4 join C6) meet C7 meet (T4 join C3));
\end{verbatim}
\noindent and the singularities are no worse than stated:
\begin{verbatim}
[ResolutionGraph(T4 join C3 join C6 join C7,p):p in P];
\end{verbatim}
To clarify the situation at the origin, we use:
\begin{verbatim}
TC:=TangentCone(T4 join C3 join C6 join C7,P2![0,0,1]);
DefiningEquation(TC) eq y*(x+3*y)*(x + 240/17*y)
*(x^2 + 82080/289*y^2)*(x^2 + 40585383/1587545*y^2);
\end{verbatim}
thus the singularity is ordinary.

\bibliography{ReferencesRito}

\

\

\noindent Carlos Rito\\
\\{\it Permanent address:}
\\ Universidade de Tr\'as-os-Montes e Alto Douro, UTAD
\\ Quinta de Prados
\\ 5000-801 Vila Real, Portugal
\\ www.utad.pt, crito@utad.pt\\
\\{\it Temporary address:}
\\ Departamento de Matem\' atica
\\ Faculdade de Ci\^encias da Universidade do Porto
\\ Rua do Campo Alegre 687
\\ 4169-007 Porto, Portugal
\\ www.fc.up.pt, crito@fc.up.pt

\end{document}